\documentclass[12pt]{article}
\usepackage{amsmath,amsfonts,amsthm,amscd,amssymb,graphicx}

\usepackage{amssymb}

\def\rit{{\Bbb R}}
\def\cit{{\Bbb C}}

\def\tit{{\Bbb T}}

\def\eps{\varepsilon}

\def\Remarks{{\noindent \it Remarks: }}

\def\bel{\begin{equation*} \begin{aligned}}
\def\eel{\end{aligned} \end{equation*}}

\def\beln{\begin{equation} \begin{aligned}}
\def\eeln{\end{aligned} \end{equation}}

\newtheorem{theorem}{Theorem}[section]

\newtheorem{e-proposition}[theorem]{Proposition}

\newtheorem{e-definition}[theorem]{Definition\rm}

\newtheorem{defi}[theorem]{Definition}

\newtheorem{theoreme}{Th\'eor\`eme}[section]

\newtheorem{proposition}[theoreme]{Proposition}

\def\og{\leavevmode\raise.3ex\hbox{$\scriptscriptstyle\langle\!\langle$~}}
\def\fg{\leavevmode\raise.3ex\hbox{~$\!\scriptscriptstyle\,\rangle\!\rangle$}}

\def\beq{\begin{equation}}
\def\eeq{\end{equation}}

\begin{document}

\centerline{\Large \bf Singularities of Rayleigh equation }

\bigskip

\centerline{D. Bian\footnote{School of Mathematics and Statistics, Beijing Institute of Technology, Beijing $100081$, China. Email: 
biandongfen@bit.edu.cn and emmanuelgrenier@bit.edu.cn}, 
E. Grenier$^1$}



\subsubsection*{Abstract}


The Rayleigh equation, which is the linearized Euler equations near a shear flow in vorticity formulation, is a key ingredient in the study of the long time behavior of solutions of linearized  Euler equations, in the study of the linear stability
of shear flows for Navier-Stokes equations and in particular in the construction of the so called Tollmien-Schlichting waves. It is also a key ingredient in the study of vorticity
depletion.

In this article we  locally describe the solutions of Rayleigh equation
near critical points of any order of degeneracy, and link  their values  on the boundary with their behaviors at infinity.


\section{Introduction}


Let $U_s(y)$ be a $C^\infty$ function on $\rit^+$, such that $U_s(y)$ converges exponentially fast, together
with its first two derivatives, to $0$. Let 
$$
\Omega = \rit \times \rit^+.
$$
We will also consider the case $\Omega = \rit \times [-1,+1]$. Without modifications, the results also extend to $\Omega = \tit \times \rit^+$
and $\tit \times [-1,+1]$ and to the corresponding three dimensional cases.
Let us consider  the shear flow 
$$
U(y) = (U_s(y), 0).
$$
The incompressible Euler equations linearized near $U$ are
\beq \label{Eul1}
\partial_t v + (U \cdot \nabla) v + (v \cdot \nabla) U + \nabla p = 0,
\eeq
\beq \label{Eul2}
\nabla \cdot v = 0,
\eeq
\beq \label{Eul3}
v \cdot n = 0 \quad \hbox{on} \quad \partial \Omega,
\eeq
where $n$ is a vector normal to $\partial \Omega$.  These equations must be completed with the initial value of $v_0(x,y)$ at $t = 0$.

A strategy to study the long time behavior of $v$ is to study the resolvent of these equations, namely to study solutions $v(\lambda,x,y)$ to
\beq \label{resol}
\lambda v + (U \cdot \nabla) v + (v \cdot \nabla) U + \nabla p = f,
\eeq
\beq \label{resol2}
\nabla \cdot v = 0,
\eeq
with the boundary condition (\ref{Eul3}), where $f = v_0$.

To solve (\ref{resol},\ref{resol2}), following the classical approach \cite{Schmid}, we introduce the stream function $\psi(\lambda,x,y)$ 
of $v(\lambda,x,y)$
 and take its Fourier transform in the $x$ variable, with dual Fourier variable $\alpha$,
 and its Laplace transform in time, with dual variable $\lambda = - i \alpha c$,  namely we look for solutions $v_\alpha(t,x,y)$
 of the Euler equations of the form
 $$
 v_\alpha(t,x,y) = \nabla^\perp \Bigl[ e^{i (\alpha x - c t)} \psi_\alpha(y) \Bigr] 
 $$
 where $\nabla^\perp = (\partial_y,-\partial_x)$.
Taking the curl of (\ref{resol}) and its Fourier transform in $x$, we obtain
\beq \label{Rayl}
(U_s(y) - c) (\partial_y^2 - \alpha^2) \psi_\alpha - U_s''(y) \psi_\alpha = f_\alpha
\eeq
where
$$
f_\alpha =  - {(\nabla \times f)_\alpha \over i \alpha},
$$
together with the boundary conditions
\beq \label{Rayl2}
\psi(0) = 0, \qquad \lim_{y \to + \infty} \psi(y) = 0,
\eeq
where $(\nabla \times f)_\alpha$ is the Fourier transform in $x$ of the curl of $f$.
This equation is  called the Rayleigh equation and is posed on $\rit^+$. It is a second order ordinary differential equation 
in $y$, singular if $U_s(y) - c$ vanishes, parametrized by the complex number $c$.

Rayleigh equation is thus the resolvent of the incompressible Euler equations linearized near a shear flow.
Its study is fundamental in the investigation of the linear and nonlinear stability of shear flows with respect to Euler and Navier-Stokes
equations.

One method to describe the solution $\psi_\alpha$ to (\ref{Rayl}) is to introduce the Green function $G_\alpha(x,y)$ of Rayleigh equation
defined by
$$
(U_s(y) - c) (\partial_y^2 - \alpha^2) G_\alpha - U_s''(y) G_\alpha = \delta_x.
$$
The solution of (\ref{Rayl},\ref{Rayl2}) is then directly given by
$$
\psi_\alpha(y) = \int_0^{+ \infty} G_\alpha(x,y) f_\alpha(x) \, dx .
$$
As $U_s(y)$ converges exponentially fast to $0$, there exists 
 two independent solutions $\psi_\pm(y,c)$ of the homogeneous Rayleigh equation (without source term),
 such that $\psi_\pm(y,c) \sim e^{\pm | \alpha | y}$ at infinity.
The construction of $G_\alpha(x,y)$ is then explicit and recalled in section \ref{secGreen}.

In this article, we focus on the construction of these two solutions $\psi_\pm(y,c)$ and study their smoothness in $y$ and $c$.
Let us from now on fix some $\alpha > 0$ and drop this index from the notations.
Of course, in the case $ f= 0$, (\ref{Rayl}) may be rewritten under the form 
\beq \label{Rayl3} 
\partial_y^2 \psi(y) = {U_s''(y) \over U_s(y) - c} \psi(y) + \alpha^2 \psi(y),
\eeq
which underlines the role of points $y$ where $U_s(y) = c$. Such couples $(y,c)$ will be called "critical layers". 

\begin{defi} 
Let $(y_0,c_0) \in \rit \times \cit$.
Then $(y_0,c_0)$ is called a critical point if $U_s(y_0) = c_0$.
It is called a critical point of order $n$  if $\partial_y^k U_s(y_0) = 0$ for any $1 \le k < n-1$ and
$\partial_y^n U_s(y_0) \ne 0$.
\end{defi}

Note that if $(y_0,c_0)$ is a critical point, then $c_0$, which equals $U_s(y_0)$, is a real number.
Near a critical layer, we expect $\psi_\pm$ to become singular. This paper is devoted to a detailed study of these singularities.

The study of solutions of (\ref{Rayl3}) near critical points $(y_0,c_0)$ of order $1$ is well known if $U_s(y)$ is analytic \cite{Reid,Reid2,Schmid}.
We then know that, locally, (\ref{Rayl3}) has then two independent solutions, one $\psi_1$, which is smooth
and of the form
$$
\psi_1(y,c) = (U_s(y) - c) P(y,c)
$$
where $P(y,c)$ is holomorphic in the vicinity of $(y_0,c_0)$, and a second one $\psi_2(y,c)$, which is singular when $U_s(y) = c$
and is of the form
$$
\psi_2(y,c) = Q(y,c) + (U_s(y) - c) R(y,c) \log( U_s(y) - c),
$$
where both $Q$ and $R$ are holomorphic in the vicinity of $(y_0,c_0)$.

If $(y_0,c_0)$ is a critical point of order $n \ge 2$, then the local structure of the solutions of the homogeneous Rayleigh equation
is up to now  unknown, even in the holomorphic case. A first study of the behavior of Rayleigh equation near singular points of order $2$ has been 
started in \cite{Bouchet}, but without a precise description of the corresponding singularity.

Critical points $(y_0,c_0)$ of order $2$ lead to the so called "vorticity depletion", a phenomena discovered by F. Bouchet and H. Morita 
in \cite{Bouchet} and then studied in mathematics for instance in \cite{Jia,Ionescu-Jia-CMP-2020,Ionescu,Wei1,Zhang,Wei2} to quote only a few results. 
To get optimal bounds on the vorticity depletion,
it appears necessary to have a precise description of the behavior of $\psi_\pm$ near the critical point \cite{BGDepletion}.

The aim of this article is to describe the solutions of Rayleigh equation near critical points of any order, for $C^\infty$ shear layer profiles $U_s(y)$. 

\begin{theorem}  \label{theo1}
Let $(y_0,c_0) \in \rit \times \cit$ be a critical point of order $n \ge 1$. Let $| \alpha | > 0$ be fixed.
Then for $(y,c)$ in the vicinity of $(y_0,c_0)$,
there exists two independent solutions $\psi_1(y,c)$ and $\psi_2(y,c)$  of Rayleigh equation, which are $C^\infty$ in $y$ and $c$ 
except at $(y,c) = (y_0,c_0)$ and such that 
\beq \label{bound1}
| \psi_1(y,c) | \lesssim 1, \quad
| \partial_y^2 \psi_1(y,c)  | \lesssim | U_s(y) - c |^{-1},
\eeq
and similarly for $\psi_2(y,c)$.
Moreover,

\begin{itemize}

\item if $n = 1$, $\psi_1(y,c)$ is  $C^\infty$ in $(y,c)$ even near $(y_0,c_0)$
and
\beq \label{bound2}
| \partial_y \psi_2(y,c) | \lesssim | \log( U_s(y) - c) |,
\eeq

\item if $n \ge 2$, $\psi_1(y,c)$ and $\psi_2(y,c)$ satisfy
\beq \label{bound3}
|\psi_1(y,c) | \lesssim |c| \zeta \Bigl( {y \over |c|^{1/n}} \Bigr), \qquad
|\psi_2(y,c) | \lesssim |c| \zeta \Bigl( -{y \over |c|^{1/n}} \Bigr),
\eeq
where
\beq \label{defizeta}
\zeta(z) = (1 + z)^n 1_{z \ge 0} + (1 + z)^{1 - n} 1_{z \le  0}.
\eeq
\end{itemize}

\end{theorem}

Bounds which are uniform in $\alpha$  are discussed in Section \ref{uniform}.
We note that $\psi_1$ and $\psi_2$ satisfy the following  "localization property":
if $y < 0$,
\beq \label{bound4}
\lim_{c \to 0} |c|^{2 - 1/n}| \psi_1(y) |  = 0
\eeq
and if $y > 0$,
\beq \label{bound5}
\lim_{c \to 0} |c|^{2 - 1/n}| \psi_1(y) |  \ne 0,
\eeq
and conversely for $\psi_2$. Asymptotically, $\psi_1$ and $\psi_2$ get "localized" on each side of $0$.

As shown in \cite{BGDepletion}, in the case $n = 2$, this localization property 
is  the key ingredient of the vorticity depletion phenomena \cite{Bouchet}.
Note that when $n = 2$ and $U_s(y) = (z-a)(z-b)$ where $a$ and $b$ are two real numbers, the Rayleigh equation is a particular
form of hypergeometric equation.

Note also that the difference between $n = 1$ and $n \ge 2$ may be interpreted in terms of monodromy. Let us assume that $U_s(y)$
is analytic, namely that it can be extended to part of the complex plane. 
Then, when $n = 1$, the Rayleigh equation in the form (\ref{Rayl3}) has only one singularity, at $y = y_c$, defined
by $U_s(y_c) = c$. When $y$ is close to $y_c$, $U_s''(y) / (U_s(y) - c)$ is large. However, for analytic functions, we can choose
to go into the complex plane and "make the turn" of the singularity by integrating Rayleigh equation for instance along a part of 
the circle of centre $y_c$
and radius $\sigma > 0$ (small but independent on $c$). On this part of circle, Rayleigh equation is not singular. As a consequence,
the solution of Rayleigh equation remains bounded on both sides of the singularity.

For $n \ge 2$ on the contrary, the contour of integration must pass through the $n$ complex zeroes of $U_s(y) - c = 0$. Near these
zeros, $U_s''(y) / (U_s(y) - c)$ is large  and, this time, when we go through these zeroes, the solutions $\psi_1$ and $\psi_2$ see
a large change in their magnitude, of order $|c|^{2 - 1/n}$.

\medskip

Once we have local solutions of Rayleigh equation, we can construct  two global solutions $\psi_\pm$. A key point in the study of
the spectrum of Orr-Sommerfeld equation is then to evaluate the ratio $\partial_y \psi_-(0,c) / \psi_-(0,c)$
where $\psi_-$ is the solution of Rayleigh such that $\psi_-(y,c) \sim e^{- | \alpha | y}$ when $y \to + \infty$.

\begin{theorem} \label{theo3}
For small $|\alpha|$, there exists two independent solutions of Ray\-leigh equation, denoted by $\psi_{\pm}(y,c)$,  defined
on  $\rit^+$, with unit Wronskian, which behave like $O(e^{\pm | \alpha | y})$ when $y$ goes to infinity. 
Moreover, for small $\alpha$ and $c$, 
\beq \label{psippsi}
{\psi_-'(0) \over \psi_-(0)} = - {U_s'(0) \over c} - {\alpha \over c^2}  (U_+ - c)^2
+ {\alpha^2 \over c^2} (U_+ - c)^4 \Omega_0(0,c) + O \Bigl( { \alpha^3 \over c^2} \Bigr) ,
\eeq  
where
\beq \label{defiOmega0}
\Omega_0(0,c) =  -{1 \over ( U_+ - c)^2} \int_0^{+ \infty} \Bigl[
{(U_s(z) - c)^2 \over (U_+ - c)^2} - {(U_+ - c)^2 \over (U_s(z) - c)^2} \Bigr] dz .
\eeq
\end{theorem}

A similar result in $[-1,+1]$ for even shear flow profiles and even solutions of Rayleigh equation (see Section \ref{interval}).

\medskip

Let us end this introduction by the study of the adjoint of Rayleigh operator, which is
\beq \label{adjoint}
Ray^t(\psi) = (\partial_y^2 - \alpha^2) (U_s - c) \psi - U_s'' \psi  .
\eeq
Then $Ray^t$ is conjugated to the Rayleigh operator
\beq \label{Rayop}
Ray(\psi) = (U_s - c) (\partial_y^2 - \alpha^2) \psi - U_s'' \psi
\eeq
in the sense that
\beq \label{conj}
Ray^t(\psi) = {1 \over U_s - c} Ray \Big( (U_s - c) \psi \Bigr) .
\eeq
Using this conjugation, all the previous theorems can be transposed into similar results for the adjoint of Rayleigh operator.

\medskip

The rest of this paper is organized as follows.
Section $2$ is devoted to the construction of solutions of Rayleigh equation near critical points,  when $\alpha = 0$ and  when
$\alpha \ne 0$, and to the study of uniform bounds. Section $3$ is devoted to the proof of Theorem \ref{theo3}, and, in Section $4$, we recall the construction
of the Green function for the Rayleigh equation in the half line.


\subsubsection*{Notation}


Throughout this paper,
$$
\langle t \rangle = 1 + | t | .
$$


\section{Local construction of solutions}


In this section we only consider the Rayleigh equation (\ref{Rayl3}), without its boundary conditions.
The Rayleigh equation is a second order differential equation in $y$, which depends on the complex parameter $c$.
It is singular at critical points. The aim of this section is to  construct locally
two independent solutions of (\ref{Rayl3}) and to study their smoothness near a critical point $(y_0,c_0)$.

Note that if $(y_0,c_0)$ is not a critical point, then we immediately get the existence of two solutions which are
$C^\infty$ in $(y,c)$ in the vicinity of $(y_0,c_0)$.


\subsection{The case $\alpha = 0$}


When $\alpha = 0$, the Rayleigh equation notably simplifies in
\beq \label{Rayalpha}
(U_s - c) \partial_y^2 \psi = U_s'' \psi 
\eeq
which has an explicit solution
\beq \label{sol1}
\psi_0(y,c) = U_s(y) - c.
\eeq
Note that $\psi_0(y,c)$ is $C^\infty$ in $(y,c)$.

An independent solution $\psi(y,c)$ is obtained through the method of the variation of the constant
\beq \label{sol2}
\psi(y,c) = (U_s(y) - c) \int_{A}^y {dz \over (U_s(z) - c)^2}
\eeq
where $A$ can be arbitrarily chosen. We choose $A > y_0$, close to $y_0$.
Let
$$
I(y,c) = \int_A^y {dz \over (U_s(z) - c)^2} .
$$
We note that the integrand is singular when $U_s(z) = c$, and in particular at the critical point $(y_0,c_0)$.
We now study in details the behavior of $I(y,c)$ near $(y_0,c_0)$.

Up to a change of variables, we may assume that $y_0 = 0$, and up to a change of $U_s(y)$ in $U_s(y) - U_s(y_0)$, we may assume that
$c_0 = 0$.


\subsubsection{Critical points of order $1$ \label{critic1}}


In this case, locally near $y_0$, $U_s'(y)$ does not vanish, thus we can make the change of variable $t = U_s(y)$.
Let
$$
V(t) = {1 \over U_s'(y)} = {1 \over U_s'( U_s^{-1}(t))}.
$$
We have, denoting by $C_0$ the various constants appearing during the following two integrations by parts, 
\bel
I(y,c) & = \int_{U_s(A)}^{U_s(y)} {V(t) \over (t- c)^2} \, dt
\cr & = - {V(U_s(y)) \over U_s(y) - c} +   \int_{U_s(A)}^{U_s(y)} {V'(t) \over t- c} \, dt + C_0
\cr & =-  {V(U_s(y)) \over U_s(y) - c} + V'(U_s(y)) \log(U_s(y) - c) 
\cr & \qquad -    \int_{U_s(A)}^{U_s(y)} V''(t) \log( t- c) \, dt + C_0.
\eel
Repeatedly integrating by parts, we obtain an expansion of $I(y,c)$ with an arbitrarily precision.
In particular, $\psi_2(y,c) = (U_s(y) - c) I(y,c)$ is bounded.
Using (\ref{Rayl3}), we obtain (\ref{bound1}), and integrating (\ref{bound1}), we obtain (\ref{bound2}) in the particular case $\alpha = 0$.


\subsubsection{Critical points of order $n \ge 2$}


Near $y_0$,  we  write $U_s(y)$ under the form
$$
U_s(y) = y^n [1 + y W(y)]
$$
where $W \in C^{\infty}$.
Let $y = |c|^{1/n} u$ and  $c = |c| e^{i \theta}$.
Then
$$
I(y,c) = |c|^{-2 + 1/n} \int_{A  |c|^{-1/n}}^{y |c|^{-1/n}} {du \over ( u^n  \Theta(u) - e^{i\theta})^2 },
$$
where
$$
\Theta(u) = 1 + |c|^{1/n} u W(| c|^{1/n} u) .
$$
We next make the change of variables
$$
t = u \Theta^{1/n}(u) .
$$
We have
$$
I(y,c) =  |c|^{- 2 + 1/n} \int_{A |c|^{-1/n} \Theta^{1/n}(A |c|^{-1/n})}^{y |c|^{-1/n} \Theta^{1/n}(|c|^{-1/n} y)}
{\Psi(t) \over (t^n - e^{i \theta})^2 } dt,
$$
where 
$$
\Psi(t) = \Bigl[ \Theta^{1/n}(u) + {u \over n} {\Theta'(u) \over \Theta^{1 - 1/n}(u)} \Bigr]^{-1}.
$$
We note that  $\Psi(0) = 1$ and that, for any $k \ge 1$,
\beq \label{boundPsi}
|\partial_t^k \Psi(t) | \lesssim |c|^{k/n}.
\eeq
Let $N$ be a large integer.
Let 
$$
P(t) = \sum_{j=0}^{2N-1} a_j t^j
$$
be the polynomial of degree $2N-1$ such that
$$
\partial_t^k P(\pm 1) = \partial_t^k \Psi(\pm 1)
$$
for any $0 \le k \le N-1$.
Then 
\beq \label{expandPsi}
\Psi(t) = P(t) + (t^2 - 1)^N \Psi_N(t),
\eeq
where $\Psi_N(t)$ is a smooth function. We note that 
\beq \label{ak}
a_j = O(|c|^{j/n}) .
\eeq
We now split $I(y,c)$ in  
\beq \label{decompI}
I(y,c) = I_P + I_R,
\eeq
where
$$
I_P =    |c|^{- 2 + 1/n} \int_{A |c|^{-1/n} \Theta^{1/n}(A |c|^{-1/n})}^{y  |c|^{-1/n}  \Theta^{1/n}(|c|^{-1/n} y)}
{P(t) \over (t^n - e^{i \theta})^2 } dt,
$$
and
$$
I_R =    |c|^{- 2 + 1/n} \int_{A |c|^{-1/n} \Theta^{1/n}(A |c|^{-1/n})}^{y  |c|^{-1/n}  \Theta^{1/n}(|c|^{-1/n} y)}
{(t^2 - 1)^N \over (t^n - e^{i \theta})^2 } \Psi_N(t) \,  dt.
$$
The first integral $I_P$ can be computed explicitly since the integrand is a rational fraction.
We decompose it in simple elements
\beq \label{decomposition}
{P(t) \over (t^n - e^{i \theta})^2} = \sum_{k=1}^n {\alpha_k \over (t - e^{i \theta_k})^2 }
+ {\beta_k \over t - e^{i \theta_k}} + Q(t)
\eeq
where
$$
\theta_k = {\theta \over n} + {2 i k \pi \over n}
$$
and where $Q(t)$ is a polynomial of degree $2N - 1 - 2n$.
It can be expanded under the form
$$
Q(t) = \sum_{k=1}^{2N - 1 - 2n} \gamma_k t^k
$$
where
$$
\gamma_k  = O( |c|^{2 + k/n}).
$$
Thus
\beq \label{boundQ}
|Q(t) | \lesssim |c|^2 + |c|^{(2N - 1)/n} |t|^{2N - 2n - 1}.
\eeq
Let $R(t) = t^n - e^{i \theta}$,  $t_k = e^{i \theta_k}$ and  $\delta = t - t_k$. Then
\bel
{P(t) \over R^2(t)} &= { P(t_k) + \delta P'(t_k) + \cdots \over \Bigl( \delta R'(t_k) + {\delta^2 \over 2} R''(t_k) + \cdots \Bigr)^2}
\cr & ={1 \over \delta^2 R'(t_k)^2}
\Bigl[ P(t_k) + \delta P'(t_k) + \cdots \Bigr] \Bigl[ 1 - \delta {R''(t_k) \over R'(t_k)} + \cdots \Bigr]
\cr & = {1 \over \delta^2} {P(t_k) \over R'(t_k)^2} 
+ {1 \over \delta} {1 \over R'(t_k)^3} \Bigl[ P'(t_k) R'(t_k) - P(t_k) R''(t_k) \Bigr] + \cdots.
\eel
Thus, identifying the coefficients of the Laurent series of both sides of (\ref{decomposition}), we obtain
$$
\alpha_k = {P(e^{i \theta_k}) \over n^2} e^{2 i \theta_k - 2 i \theta} 
$$
and
$$
\beta_k =  {P'(e^{i \theta_k}) \over n^2} e^{2 i \theta_k - 2 i \theta}  - {n-1 \over n^2}P(e^{i \theta_k}) e^{i \theta_k - 2 i \theta}.
$$
Let $J(t)$ be a primitive of $P(t) (t^n - e^{i \theta})^{-2}$. Then we can choose
\beq \label{J}
J(t) = -  \sum_{k=1}^n {\alpha_k \over t - e^{i \theta_k}} + \sum_{k=1}^n \beta_k \log(t - e^{i \theta_k}) + Q_1(t)
\eeq
where $Q_1$ is a primitive of $Q$. 

Let us now study the behavior of $J(t)$ for large $t$.
We recall that
\beq \label{sumroot}
\sum_{k=1}^n (e^{2 i k \pi \over n})^p = n 1_{n \, | \, p}
\eeq
where $n \, | \, p$ means that $p$ is a multiple of $n$.
Thus we have
\bel 
\sum_{k=1}^n {\alpha_k \over t - e^{i \theta_k}} 
&= \sum_{p \ge 0}  {1 \over t^{p+1}} \sum_{k=1}^n \alpha_k e^{i p \theta_k}
\cr & = \sum_{p \ge 0} {1 \over n^2 t^{p+1} } \sum_{k=1}^n P(e^{i \theta_k}) e^{ i (p + 2) \theta_k - 2 i \theta}
\cr & =  \sum_{p \ge 0} {1 \over n^2 t^{p+1} } \sum_{k=1}^n \sum_{j=0}^{2N - 1} a_j  e^{ i (2 + p + j) \theta_k - 2 i \theta}
\cr & = {1 \over n} \sum_{j = 0}^{2N-1} a_j \sum_{p \ge 0, \, 2 + p + j = l n ,\,  l \ge 1} {e^{i (l-2) \theta} \over t^{p+1}}.
\eel
Let us detail this sum. 
The terms involving $l = 1$ equal
$$
P_1 = \sum_{j = 0}^{n-2} {a_{j} \over n t^{n - j - 1}} e^{- i \theta}.
$$
The terms involving $l = 2$ are of order 
$$
 O\Bigl({1 \over t^{2 n - 1}} \Bigr) + O \Bigl ({ |c|^{1/n}  \over t^{ 2n - 2}} \Bigr) + \cdots + O \Bigl({ |c|^{2 - 2/n} \over t} \Bigr) 
 = O \Bigl( {1 \over t^{2 n - 1}} \Bigr) + O \Bigl( {|c|^{2 - 2/n}  \over t} \Bigr). 
$$
The terms involving $l \ge 3$ are smaller, thus
\beq \label{firstsum}
\sum_{k=1}^n {\alpha_k \over t - e^{i \theta_k}}  = P_1 +  O(t^{-2 n +1}) + O(|c|^{2 - 2/n} t^{-1}).
\eeq
Moreover, when $t > 0$, we choose the usual determination of the logarithm, and
\bel
\sum_{k=1}^n \beta_k \log(t - e^{i \theta_k})
& =  \sum_{k=1}^n \beta_k \log t+ \sum_{k=1}^n \beta_k \log(1 - t^{-1} e^{i \theta_k}).
\eel
We note that
\bel
  \sum_{k=1}^n  \beta_ke^{i p \theta_k}  & =
  \sum_{j=0}^{2N-1}  \sum_{k=0}^n
  a_j \Bigl[ {j - n+1 \over n^2} \Bigr] e^{i (p + j + 1) \theta_k - 2 i \theta}
 \cr & = \sum_{j = 0}^{2N-1} {(l-1) n - p \over n} a_j e^{i (l - 2) \theta} 1_{n \, | \, p + j + 1 },
 \eel
 where in this formula, $l$ is defined by $p + j + 1 = l n$.
This leads to
\bel
 \sum_{k=1}^n \beta_k    =  \sum_{l \ge 1} (l-1) a_{l n -1} e^{i (l - 2) \theta} = O(|c|^{2 - 1/n}) .
 \eel
Moreover,
\bel
\sum_{k=1}^n \beta_k \log(1 - t^{-1} e^{i \theta_k}) 
& = - \sum_{p \ge 1} {1 \over p t^p}  \sum_{k=1}^n  \beta_ke^{i p \theta_k}
\cr & = - \sum_{j = 0}^{2N - 1} a_j \sum_{p \ge 1} {(l-1) n - p \over  p t^p} e^{i (l - 2) \theta}   1_{n | p + j + 1}.
\eel
The terms corresponding to $l = 1$ are
$$
\sum_{j = 0}^{n-2} a_j {e^{- i \theta} \over t^{n - j - 1}} =  P_1.
$$
As previously, the terms involving $l \ge 2$ are of order
$$
O\Bigl( {1 \over t^{2 n - 1}} \Bigr) + \cdots + \Bigl( {|c|^{2 - 2 / n} \over t} \Bigr).
$$
When $t < 0$, if $\Im e^{i \theta_k} < 0$,
$$
\log( t - e^{i \theta_k}) = \log( (-t) + e^{i \theta_k}) - i \pi,
$$
else
$$
\log( t - e^{i \theta_k}) = \log( (-t) + e^{i \theta_k}) + i \pi .
$$
Thus
\beq \label{minusinfty}
\sum_{k=1}^n \beta_k \log(t - e^{i \theta_k})
= \Gamma + \sum_{k=1}^n \beta_k \log( (-t) + e^{i \theta_k})
\eeq
where
$$
\Gamma = - i \pi \sum_{k = 1}^n \beta_k sign(\Im e^{i \theta_k}) \ne 0.
$$
The second term of (\ref{minusinfty}) can be bounded as previously.

It remains to compute the contribution of $Q_1$ to $J$. We note that $Q$ is a polynomial of degree at most $2N - 1 - n$.
Choosing a primitive $Q_1$ of $Q$ which vanishes at $t = 0$ and using (\ref{boundQ}), we obtain
$$
| Q_1(t) | \lesssim |c|^2 |t| + |c|^{(2N - 1)/n} t^{2 N - 2n}.
$$
We have
\beq \label{IP1}
I_P =   |c|^{- 2 + 1/n} \Bigl[ J \Bigl( y |c|^{-1/n}\Theta^{1/n}(| c|^{-1/n} y ) \Bigr) -   J \Bigl( A |c|^{-1/n} \Theta^{1/n}(| c|^{-1/n} A ) \Bigr) \Bigr].
\eeq
However, the integral (\ref{sol2}) defining $\psi$ is defined up to a constant. We thus integrate the second term of (\ref{IP1}) in this constant,
and omit the second term in  (\ref{IP1}). With this new definition of $I(y,c)$, $\psi$ is still a solution of Rayleigh equation.

We also have to take care of the term in $\log t$. By adding a constant, we replace
$$
\sum_{k=1}^n \beta_k \log t
$$
by 
$$
\sum_{k = 1}^n \beta_k \log \Bigl( {t \over |c|^{1/2}} \Bigr).
$$
With this new definition of $I_P$, for $t > 0$, we have
\bel
| (U_s - c) I_P |  & \lesssim |U_s(y) - c | | c |^{-2 + 1/n} 
\Bigl[ {1\over \langle t \rangle^{2n-1}} + {|c|^{2 - 2/n} \over \langle t \rangle } + |c|^{2 - 1/n} \log \Bigl( 2 + { |t| \over |c|^{1/n}}  \Bigr) \Bigr]
\cr & \quad + | U_s(y) - c | |c|^{-2 + 1/n} \Bigl[ |c|^2 |t| + |c|^{(2N - 1)/n} |t|^{2N - 2n} \Bigr]
\cr & \lesssim {1 \over |c|^{1 - 1/n} \langle t \rangle^{n-1}} + |c|^{1-1/n} \langle t \rangle^{ n - 1} 
+ |c| \langle t \rangle^n \log \Bigl( 2 +  {|t| \over |c|^{1/n}} \Bigr) \rangle
\cr & \quad + |c| |t|^n \Bigl[ |c|^{1/n} |t| + |c|^{(2N - 2n)/n} |t|^{2N - n} 
\cr & \lesssim {1 \over (|y| + |c|^{1/n})^{n-1}} + (|y| + |c|^{1/n})^{n-1} + (|y| +  |c|^{1/n})^n \log (2 + |y| )
\cr & \quad + (|c|^{1/n} + |y|)^{n+1} + (|c|^{1/n} + |y|)^{2N}
\cr & \lesssim {1 \over (|y| + |c|^{1/n})^{n-1}} + (|y| + |c|^{1/n})
\cr & \lesssim {1 \over (|y| + |c|^{1/n})^{n-1}} ,
\eel
provided $|y|$ and $|c|$ are small enough,
where we have used that
$$
U_s(y) - c = y^n \Bigl[1 + y W(y) \Bigr] - c
= c u^n \Bigl[1 + y W(y) \Bigr] - c
$$
which implies that
$$
| U_s(y) -c | \lesssim |c| t^n + |c| \lesssim  |c| \langle t \rangle^n.
$$
For $t < 0$, we have
\bel
|  (U_s(y) - c) I_P | &\lesssim | U_s(y) - c| |c |^{-2 + 1/n}
 \lesssim {\langle t \rangle^{n} \over |c|^{1 - 1/n}}
 \lesssim {(|y| + |c|^{1/n})^n \over |c|^{2 - 1/n}} 
\eel
since the other contributions are negligible with respect to this leading term.

We now turn to $I_R$. The main observation is that
\beq \label{claim}
\Bigl| { (t^2 - 1)^n \over (t^n - e^{i \theta})^2} \Bigr| \lesssim 1.
\eeq
Namely, if $n$ is even,
\bel 
\Bigl| {(t^2 - 1)^n \over (t^n - e^{i \theta})^2 } \Bigr|
= \Pi_{k = -n/2}^{n/2-1} \Bigl| {t - 1\over t - e^{i \theta_k}} \Bigr|^2
\Pi_{k = n/2}^{3n/2-1} \Bigl| {t + 1\over t - e^{i \theta_k}} \Bigr|^2.
\eel
Let $- n/2 \le k \le n/2 -1$. Then $\arg(\theta_k) \in [-\pi/2,\pi/2]$. We want to prove
$$
\Bigl| {t - 1\over t - e^{i \theta_k}} \Bigr| \lesssim 1.
$$
It is true if $t \notin [1/2,3/2]$ since then the denominator is bounded away from $0$.
If $t \in [1/2,1]$, the circle of centre $t$ and radius $1 -t$ is inside the circle of centre $0$ and radius $1$.
Thus, 
\beq \label{elem}
|t - e^{i \theta_k} | \ge |t - 1|.
\eeq
 If $t \in [1, 3/2]$, then the circle of centre $t$ and radius $1 - t$ is outside
the circle of centre $0$ and radius $1$, thus (\ref{elem}) is also true.
Combining these various cases, we obtain (\ref{claim}). The case where $n$ is odd is similar.

Note that (\ref{claim}) implies that the integrand of $I_R$ is not singular, thus $I_R$ is bounded. Derivatives of $I_R$, up to 
the order $N - n$ are also bounded.
Let
$$
\psi_2(y,c) = c^{3 - 1/n} (I_P + I_R) .
$$
Then $\psi_2$ is bounded and satisfies (\ref{bound5}). Using Rayleigh equation we further obtain (\ref{bound1}).

Let 
$$
\langle y \rangle_c = |c|^{1/n} + |y|,
$$
and let us define $\zeta(y,c)$ by
\beq \label{zeta1}
\zeta(y,c) = \left\{ \begin{array}{cc} \langle y \rangle_c^n
&\quad \hbox{if} \quad y > 0 ,
\cr
 |c|^{2 - 1/n} \langle y \rangle_c^{-n+1} 
 &\quad \hbox{if} \quad y < 0.
 \end{array} \right.
\eeq
Then
\beq \label{estimpsi2}
| \psi_2(y,c) | \lesssim \zeta(-y,c).
\eeq
The construction of $\psi_1$ is similar, except that we use the determination of the logarithm
which is defined on $\cit - \rit^+$ instead of the usual one. We then have
\beq \label{estimpsi1}
| \psi_1(y,c) | \lesssim \zeta(y,c).
\eeq


\subsection{The case $\alpha \ne 0$ \label{nonzero}}


In view of (\ref{Rayl2}), near the critical point, $\alpha^2$ is negligible with respect to $U_s'' / (U_s - c)$ and can be treated
as a perturbation.
 Let $G_0$ be the Green function of 
$$
\partial_y^2 \phi =   {U_s'' \over U_s - c} \phi + \delta_x.
$$
Then $G_0$ is explicitly given by 
 \beq \label{defiG1}
G_0(x,y) = W^{-1} \Bigl\{
\begin{array}{c} \psi_1(x) \psi_2(y) \quad \hbox{if} \quad x < y, \cr
\psi_1(y) \psi_2(x) \quad \hbox{if} \quad x > y,
\end{array} 
\eeq
where $W$ is the Wronskian of $\psi_1$ and $\psi_2$ is of order $O(|c|^{2 - 1/n})$.

Let $\sigma > 0$. For any bounded function $f$, we define
\beq \label{defiGreen1}
 {\cal G} f(y) =  \int_{-\sigma}^\sigma G_0(x,y) f(x) \, dx.
\eeq
Let us first bound ${\cal G}$. 
We define the norm $\| \phi \|_\zeta$ by
$$
\| \phi \|_\zeta = \sup_{-\sigma \le x \le \sigma} | \zeta(x,c)|^{-1}  \, | \phi(x)| .
$$
We have
$$
{\cal G}f(y)  
  =  \psi_2(y)  \int_{-\sigma}^y  { \psi_1(x)  \over W} f(x)  \, dx 
+  { \psi_1(y) \over W}  \int_y^{+\sigma}  \psi_2(x)   f(x)  \, dx.
$$
Let us bound the first term. For $y < 0$, we have
\bel
\Bigl|  \psi_2(y)  \int_{-\sigma}^y  { \psi_1(x)  \over W} f(x)  \, dx \Bigr|
& \lesssim {\langle y \rangle_c^n \over |c|^{2 - 1/n}} \int_{-\sigma}^y 
\Bigl[ {|c|^{2 - 1/n} \over \langle x \rangle_c^{n - 1}}   \Bigr]^2 \, dx
\cr & \lesssim \langle y \rangle_c^n {|c|^{2 - 1/n} \over \langle y \rangle_c^{2 n - 3}} 
  \lesssim  \sigma^2  {|c|^{2 - 1/n} \over \langle y \rangle_c^{ n - 1}} 
\cr & \lesssim  \sigma^2 \zeta(y)
\eel
and for $ y > 0$, the same integral is bounded by
\bel
{1 \over |c|^{2 - 1/n}} {|c|^{2 - 1/n} \over \langle y \rangle_c^{n-1}}
& \Bigl[ \int_{-\sigma}^0 { |c|^{4 - 2/n} \over \langle x \rangle_c^{2 n - 2} } \, dx
+ \int_0^y \langle x \rangle_c^{2n} \, dx \Bigr]
\cr &  \lesssim 
{1 \over \langle y \rangle_c^{n-1}} \Bigl[ |c|^{4 - 2/n} + \langle y \rangle_c^{2 n + 1} \Bigr]
\lesssim   \sigma^2  \langle y \rangle_c^n \lesssim    \sigma^2 \zeta(y).
\eel
We now turn to the second integral. We have, for $y > 0$,
\bel
\Bigl| { \psi_1(y) \over W}  \int_y^{+\sigma}  \psi_2(x)   f(x)  \, dx \Bigr|
 & \lesssim {\langle y \rangle_c^n \over |c|^{2 - 1/n}} \int_y^{+ \sigma}
 {|c|^{2 - 1/n} \over \langle x \rangle_c^{n - 1}} \langle x \rangle_c^n  \, dx
 \cr \lesssim   \sigma^2 \langle y \rangle_c^n \lesssim    \sigma^2 \zeta(y),
 \eel
 and for $y < 0$, this integral is bounded by
\bel
{1 \over |c|^{2 - 1/n}} {|c|^{2 - 1/n} \over \langle y \rangle_c^{n-1}}
& \Bigl[ \int_0^\sigma {|c|^{2 - 1/n} \over \langle x \rangle_c^{n-1}} \langle x \rangle_c^n \, dx
+ \int_{-y}^0 \langle x \rangle_c^n {|c|^{2 - 1/n} \over \langle x \rangle_c^{n-1}} \, dx \Bigr]
\cr & \lesssim    \sigma^2 {|c|^{2 - 1/n} \over \langle y \rangle_c^{n-1}} \lesssim    \sigma^2 \zeta(y).
\eel
Thus, for $|c|$ small enough, we have
\beq \label{eestG}
  \| {\cal G} f(y) \|_\zeta \lesssim    \sigma^2   \| f \|_\zeta.
\eeq
We  define a first solution $\psi_{1,\alpha}$ of Rayleigh equation with $\alpha \ne 0$  through the iterative scheme
\beq \label{iter00}
 \phi_{n+1} = \psi_1 + \alpha^2 {\cal G}  \phi_n,
\eeq
starting with $\phi_0 = \psi_1$.
Then, provided $\alpha^2 \sigma^2$ is small enough,  $\phi_n$ converges in the norm $\| \cdot \|_\zeta$ to a solution
 $\psi_{1,\alpha}$ of Rayleigh equation.
 Symmetrically, we construct another independent solution $\psi_{2,\alpha}$ which satisfies symmetric bounds.
 We note that
 $$
| \psi_{1,\alpha}(y,c) | \lesssim \zeta(y), \qquad
| \psi_{2,\alpha}(y,c) | \lesssim \zeta(-y).
$$
In particular, $\psi_{1,\alpha}$ and $\psi_{2,\alpha}$ are bounded. This ends the proof of Theorem \ref{theo1}.


\subsection{Uniform estimates \label{uniform}}


We now discuss how to obtain  estimates which are uniform in $\alpha$.
Two scales appear when $|\alpha|$ is large, first $|c|^{1/n}$, which is the modulus of the singularities and second $| \alpha |^{-1}$,
which is the "scale of evolution" of $\partial_y^2 - \alpha^2$.
Hence, three cases appear, depending on whether $|c|^{1/n} \ll | \alpha |^{-1}$, $|c|^{1/n} \approx| \alpha |^{-1}$ or $|c|^{1/n} \gg |\alpha|^{-1}$.

\begin{proposition}
Let
$$
K(y) = \alpha^2 + {U_s''(y) \over U_s(y) - c }.
$$
 Then, there exist $\eps > 0$, $\sigma > 0$, $\sigma_0 > 0$ and $\vartheta > 0$ such that, uniformly  in $\alpha \in \rit$ and uniformly in $|c| \le \eps$, there exist two solutions 
$\psi_+(y)$ and $\psi_-(y)$ of Rayleigh equation which satisfy

\begin{itemize}

\item if $| \alpha c^{1/n} | \le \sigma$,
 
\begin{itemize}
\item for $ y_0 = \sigma_0 |\alpha|^{-1} \le y \le \vartheta$,
\beq \label{uniformest1}
|\psi_-(y)|  \lesssim 
\exp \Bigl( - \int_{y_0}^y \Re \sqrt{K(z)} \, dz   \Bigr),
\eeq
\item for $-y_0 \le y \le y_0$,
\beq \label{uniformest2}
| \psi_- (y) | \lesssim {\zeta(- |c|^{-1/n} y) \over \zeta(- |c|^{-1/n} y_0)},
\eeq
\item for $- \vartheta  \le y \le  - y_0$,
\beq \label{uniformest3}
|\psi_-(y)|  \lesssim  {\zeta(|c|^{-1/n} y_0) \over  \zeta(- |c|^{-1/n} y_0)} 
\exp \Bigl( - \int_{-y_0}^y \Re \sqrt{K(z)} \, dz   \Bigr),
\eeq
\end{itemize}

\item if $| \alpha c^{1/n} | \ge \sigma$ then, for $- \vartheta \le  y  \le \vartheta$,
\beq \label{uniformest4}
|\psi_-(y)|  \lesssim  \exp \Bigl( - \int_0^y \Re \sqrt{K(z)} \, dz   \Bigr),
\eeq
\end{itemize}

and symmetrically for $\psi_+(y,c)$.

\end{proposition}

\Remarks
When $| \alpha|^{-1} \gg |c|^{1/n}$, between $y = y_0$ and $y = - y_0$,  there is an "amplification factor"
$\zeta(|c|^{-1/n} y_0) / \zeta(- |c|^{-1/n} y_0)$ which is of order
$$
 \Bigl( 1 + {1 \over |\alpha c^{1/n} |} \Bigr)^{2n - 1}.
$$
This amplification factor decreases when $|\alpha|$ increases and is of order $1$ when $|\alpha|^{-1}$ is of order $|c|^{1/n}$.
This makes  the transition between (\ref{uniformest1},\ref{uniformest2},\ref{uniformest3}) and (\ref{uniformest4}). 
For $| \alpha | \gg |c|^{-1/n}$, there is no amplification, as in the case of degenerate points of order $1$.
This makes the transition between a degenerate point of order $\ge 2$ and two "separated"  degenerate points of order $1$.

To obtain estimates which are uniform in $\alpha$, we use the classical WKBJ method.
Away from the singularities, we look for solutions of the form $e^{\theta(y)}$
and describe in detail $\theta(y)$. Near the singularities, this description is no longer valid. We are in a situation similar to that of "turning points",
except that here we have a singularity of $K(y)$ and not a zero. We thus have to describe the solutions near singularities, and match them with $e^{\theta(y)}$.

The WKBJ method is classical, and we will only sketch the proof of this proposition.

\begin{proof}

Let us first apply the WKBJ method to the Rayleigh equation.
We look for two solutions $\psi_\pm(y)$ under the form 
$\psi_\pm(y) = e^{\theta_\pm(y)}$, which leads to
\beq \label{WKBJ1}
\theta_\pm'^2 + \theta_\pm'' = K(y).
\eeq
The analysis of (\ref{WKBJ1}) is standard and we will just state the result. Provided $K(y)$ is a small perturbation of $\alpha^2$, and more precisely, provided
\beq \label{condWKBJ}
| K' | \ll | K |^{3/2},
\eeq
we have
\beq \label{WK}
\psi_\pm(y) = {K^{1/4}(y_0) \over K^{1/4}(y)} \exp \Bigl( \pm \int_{y_0}^y \sqrt{K(z)} \, dz [1 + O(|\alpha|^{-1})]  \Bigr),
\eeq
where $y_0$ is arbitrary.
The first derivative $\partial_y \psi_\pm(y)$ has a similar expression.
Note that if $c \ne 0$, all the zeros of $U_s(y) - c$ are simple, thus the integral appearing in (\ref{WK}) is well-defined.
Moreover, $\psi_\pm(y)$ is normalized by $\psi_\pm(y_0) = 1$.
We note that $K^{1/4}(y_0)$ and $K^{1/4}(y)$ are both of order $|\alpha|^{1/2}$ if $K(y)$ is a small perturbation of $| \alpha |^2$.

Let us discuss  (\ref{condWKBJ}) in the case $U_s(y) = y^n$ in order to simplify the presentation, the general case being similar.
Then
$$
| K'(y) | \lesssim {|y|^{n-3} \over |y^n - c|} + {|y|^{2 n - 3} \over |y^n - 1|^2} .
$$
Let us define $z$ by $z = |c|^{-1/n} y$. Then (\ref{condWKBJ}) is satisfied if
\beq \label{condWKBJ2}
{|z|^{n-3} \over |z^n - 1|} + {|z|^{2 n - 3} \over |z^n - 1|^2} \lesssim (\alpha |c|^{1/n})^3.
\eeq
If $|\alpha c^{1/n}|$ is large, then (\ref{condWKBJ2}) is satisfied provided $z - 1  \gg |\alpha c^{1/n}|^{-3/2}$, namely
provided 
$$
y - |c|^{1/n} \gg {1 \over | \alpha | | \alpha c^{1/n}|^{1/2}},
$$
If $|\alpha c|^{1/n}$ is small,  (\ref{condWKBJ}) is satisfied provided $z \gg  | \alpha c^{1/n}|^{-1}$, namely provided $y \gg | \alpha |^{-1}$.

In the sequel we will distinguish the cases $|\alpha c^{1/n}| \le \sigma$, $\sigma \le | \alpha c^{1/n} | \le \sigma'$ and $\sigma' \le | \alpha c^{1/n}|$,
where $\sigma$ (small) and $\sigma'$ (large) will be chosen later.


\subsubsection*{The case $| \alpha c^{1/n}|  \le \sigma$}


Let $\sigma$ be small enough, to be fixed below.
The construction of the Section \ref{nonzero} is valid on $[-y_0, y_0]$ where $y_0 =  \sigma_0 |\alpha|^{-1}$ provided $\sigma_0$ is small enough
and gives two solutions $\psi_1$ and $\psi_2$. 
We note that $|c|^{1/n}$ is in this interval provided $\sigma \le \sigma_0$, which we assume from now on.

Moreover, the WKBJ method can be applied provided $|y \pm | c |^{1/n}|$ is large enough with respect to $|\alpha|^{-1}$.
For $y > y_1 = |c|^{1/n} + \sigma_1 |\alpha|^{-1}$, where $\sigma_1$ is large enough, two solutions $\psi_\pm(y)$ are given by (\ref{WK}).
We note that, at $y = y_1$, $\partial_y \psi_\pm(y_1) / \psi_\pm(y_1)$ is of order $| \alpha |$.

Now on $y_0 \le y \le y_1$, we rescale the Rayleigh equation by a factor $|\alpha|$
and study $\phi(z) = \psi( | \alpha |^{-1} z)$.
The coefficients of the rescaled equation on $\phi(z)$ are bounded, thus the corresponding resolvent is bounded uniformly in $|\alpha|$.
As a consequence,  $\psi_\pm(y_0)$ and $| \alpha |^{-1} \partial_y \psi_\pm(y_0)$ 
are bounded uniformly in terms of $\psi_\pm(y_1)$ and $| \alpha |^{-1} \partial_y \psi_\pm(y_1)$. 
As  $K(y)$ is of order $\alpha^2$ when $y \ge y_0$,  we have
\beq \label{estima1}
|\psi_\pm(y)|  \lesssim {| \alpha |^{1/2} \over K^{1/4}(y)} \exp \Bigl( \pm \int_{y_0}^y \Re \sqrt{K(z)} \, dz   \Bigr)
\eeq
when $y \ge y_0$. In particular, $\psi_\pm(y_0)$ is of order $O(1)$, and $\partial_y \psi_\pm(y_0)$ is of order $O(|\alpha |)$.

We now extend $\psi_\pm$ by a linear combination of $\psi_1$ and $\psi_2$ between $-y_0$ and $y_0$. 
We write that, for some coefficients $a$ and $b$, 
$$
\psi_- = a \psi_1 + b\psi_2
$$ 
and that
$$
a \psi_1(y_0) + b \psi_2(y_0) = \psi_-(y_0),
$$
$$
a \partial_y \psi_1(y_0) + b \partial_y \psi_2(y_0) = \partial_y \psi_-(y_0).
$$
We choose $\psi_1$ and $\psi_2$ such that
$$
| \psi_1(y) | \le Z^n 1_{y >0} + Z^{1-n} 1_{y < 0}
$$
where
$$
Z = 1 + {y \over |c|^{1/n}} ,
$$
and symmetrically for $\psi_2$.
Let 
$$
Z_0 = 1 + {y_0 \over   |c|^{1/n}} .
$$
We have
$$
\left( \begin{array}{c} a \cr b \end{array} \right) 
= {1 \over W} \left( \begin{array}{cc} 
\partial_y \psi_2(y_0) & - \psi_2(y_0) \cr
- \partial_y \psi_1(y_0) & \psi_1(y_0) \end{array} \right) 
\left(  \begin{array}{c} \psi_-(y_0) \cr \partial_y \psi_-(y_0) \end{array} \right) 
 $$
where the  Wronskian $W$ of $\psi_1$ and $\psi_2$ is of order $O(|c|^{-1/n})$.  Thus
$$
a = {1 \over W} \Bigl( \partial_y \psi_2(y_0)  \psi_-(y_0) -  \psi_2(y_0) \partial_y \psi_-(y_0) \Bigr)
$$
and
$$
b = {1 \over W} \Bigl( - \partial_y \psi_1(y_0)  \psi_-(y_0) + \psi_1(y_0) \partial_y \psi_-(y_0) \Bigr).
$$
We note that $\partial_y \psi_- (y_0)/ \psi_-(y_0)$, $\partial_y \psi_1 (y_0)/ \psi_1(y_0)$ and $\partial_y \psi_2(y_0)/ \psi_2(y_0)$ are all of order $|\alpha|$, thus,
$$
a = O \Bigl( |\alpha c^{1/n} | Z_0^{1-n} \psi_-(y_0) \Bigr), \qquad b = O \Bigl( |\alpha c^{1/n}| Z_0^n \psi_-(y_0) \Bigr).
$$
As a consequence, $a \psi_1(y_0)$ and $b \psi_2(y_0)$ are both of order $\psi_-(y_1)$.

For $-y_0 \le y \le y_0$, the solution is given by $\psi_- = a \psi_1 + b \psi_2$. In particular, 
$$
\psi_-(-y_0) = O \Bigl( | \alpha c^{1/n} | Z_0^{2 n} \psi_-(y_0) \Bigr)
$$ 
and $\partial_y \psi_-(-y_0)$ is larger by a factor $|\alpha|$. Note that $| \alpha c^{1/n} | Z_0$ is bounded.
 More generally, we obtain
\beq \label{estima2}
| \psi_- (y) | \lesssim {\zeta(- |c|^{-1/n} y) \over \zeta( - |c|^{-1/n}  y_0)}
\eeq
when $-y_0 \le y \le y_0$.

We  extend $\psi_-$ for $-y_1 \le y \le -y_0$ using again the rescaled Rayleigh equation  and we  match it with a solution given by the
WKBJ method.
We obtain, for $- \vartheta \le y \le - y_0$,
\beq \label{estim3}
|\psi_-(y)|  \lesssim  {\zeta(- |c|^{-1/n}  y) \over \zeta(- |c|^{-1/n}  y_0)} {| \alpha |^{1/2} \over K^{1/4}(y)} \exp \Bigl( - \int_{-y_0}^y \Re \sqrt{K(z)} \, dz   \Bigr).
\eeq


\subsubsection*{The case $\sigma \le | \alpha c^{1/n} | \le \sigma'$}


For $|y| \ge y_1 =  |c|^{1/n} + \sigma_1 | \alpha |^{-1}$ with $\sigma_1$ large enough, we use the WKBJ method. For $|y| \le y_1$, we rescale the Rayleigh equation by introducing $z = | \alpha |^{-1} y$, 
namely we introduce $\phi(z) = \psi(| \alpha| z)$. We have
\beq \label{RR}
\partial_z^2 \phi = \phi + {U_s''(| \alpha | z) \over \alpha^2 [U_s( | \alpha | z) - c]} \phi.
\eeq
This equation is singular near $z_\pm = \pm |\alpha c^{1/n}|$, which are of order $O(1)$. These singularities are simple, thus, using the construction of  Section  \ref{nonzero},
we can construct solutions on $[z_\pm - \sigma_0 |\alpha|^{-1}, z_\pm + \sigma_0 | \alpha |^{-1}]$, provided $\sigma_0$ is small enough.
The resolvent on these intervals are moreover bounded.
Outside these intervals, the coefficients of (\ref{RR}) are bounded, thus the resolvent is also bounded.
As a consequence, the resolvent of this equation between $\pm y_1$ is bounded, uniformly in $\alpha$.

We then match the solutions constructed for $|y | \ge y_1$ and $|y| \le y_1$, which gives (\ref{uniformest4}).


\subsubsection*{The case $| \alpha c^{1/n} | \ge \sigma'$}


Let $\sigma' \gg 1$, which will be chosen below.
In this case, $K'$ is negligible with respect to $K^{3/2}$ provided 
$$
\Bigl| y - |c|^{1/n} \Bigr| \ge \delta_1 = {\sigma_1 \over | \alpha  | \, | \alpha c^{1/n}|^{1/2}},
$$
where $\sigma_1$ is large enough.
We rescale the Rayleigh equation by introducing $z = |c|^{-1/n} y$. Let $\phi(z) = \psi(|c|^{1/n} z)$, then $\phi(z)$ satisfies
\beq \label{R2}
\partial_z^2 \phi = (\alpha |c|^{1/n})^2 \phi +  U_s(|c|^{1/n} z) \phi.
\eeq
The equation (\ref{R2}) can be solved using the WKBJ method provided $|z \pm 1 | \ge \sigma_1 | \alpha c^{1/n}|^{-3/2} $ and using the construction of  Sections \ref{critic1} and \ref{nonzero}
for $| z \pm 1 | \le \sigma_0 |\alpha c^{1/n}|^{-1}$. These two regions overlap provided $\sigma'$ is large enough.
We  thus match together the various asymptotics and obtain (\ref{uniformest4}).
\end{proof}


\section{Global construction of solutions}



\subsection{Construction in the half line} 


We observe that, as $y$ goes to $+ \infty$, the only possible asymptotic behaviours of the solutions are $e^{+ \alpha y}$ and
$e^{- \alpha y}$, assuming $\alpha > 0$ to fix the ideas.
Classical results on holomorphic ordinary differential equations provide the existence of two solutions $\psi_{+,\alpha}(y,c)$ and $\psi_{-,\alpha}(y,c)$
defined for $y \ge y_1$ with $y_1$ large enough, and  such that
$$
\psi_{\pm,\alpha}(y,c) \sim e^{- \alpha y}
$$
as $y \to + \infty$.
We then extend $\psi_{\pm,\alpha}(y,c)$ towards $y = 0$ by matching together the constructions of the previous paragraphs.

We now turn to the proof of Theorem \ref{theo3}, which makes the link between the behavior at $+ \infty$ and at $y = 0$ when $\alpha$ is small.
 We introduce, following J.W. Miles \cite{Miles},
\beq \label{definitionOmega}
\Omega(y) = {\psi \over (U_s - c) [ U_s' \psi - (U_s-c) \psi'] }
\eeq
and note that, after some calculations, $\Omega(y)$ satisfies the ordinary differential equation
\beq \label{equationOmega}
\Omega' = \alpha^2 Y \Omega^2 - Y^{-1}
\eeq
where 
$$
Y(y) = (U_s(y) - c)^2.
$$  
We note that the source term $Y^{-1}(y)$ is singular when $U_s(y) = c$.
By definition of $\Omega(y)$, for  $\psi = \psi_{-,\alpha}$, 
the corresponding function $\Omega(y)$ satisfies 
$$
\lim_{y \to + \infty} \Omega(y) = {1 \over \alpha ( U_+ - c)^2} .
$$
Using (\ref{equationOmega}), we expand $\Omega(y)$ in
$$
\Omega(y) = {1 \over \alpha (U_+ - c)^2} + \theta(y)
$$ 
where $\theta$ satisfies
\beq \label{equtheta}
\theta' = \Bigl[ {(U_s - c)^2 \over (U_+ - c)^4} - {1 \over (U_s - c)^2} \Bigr]
+ 2 \alpha {(U_s - c)^2 \over (U_+ - c)^2} \theta + \alpha^2 (U_s - c)^2 \theta^2 .
\eeq
Let
\beq \label{defiOmega0z}
\Omega_0(y,c) = -{1 \over ( U_+ - c)^2} \int_y^{+ \infty} \Bigl[
{(U_s(z) - c)^2 \over (U_+ - c)^2} - {(U_+ - c)^2 \over (U_s(z) - c)^2} \Bigr] dz . 
\eeq
Note that the integral defining in $\Omega_0(y,c)$ is well defined at $+\infty$ since $U_s(z)$ converges exponentially fast to $U_+$ at infinity.

The second term in the integral defining $\Omega_0(y,c)$  is singular if $U_s(y) = c$. 
More precisely, near a critical point $y_c$,
\beq \label{exx}
{1 \over (U_s(z) - c)^2} = {1 \over U_c'^2 (z -y_c)^2} - {U_c'' \over U_c'^3} {1 \over z - y_c} + \cdots.
\eeq
Thus $\Omega_0(y,c)$ behaves like $(y - y_c)^{-1}$ near $y_c$. As a consequence, $(U_s - c) \Omega_0(y,c)$ is bounded.
Using (\ref{equtheta}) we further expand $\Omega(y)$ in 
\beq \label{expansionOmega}
\Omega(y) = {1 \over \alpha (U_+ - c)^2} + \Omega_0(y,c) + \theta_1(y).
\eeq
Inserting in (\ref{equationOmega}), we obtain after a few computations that
$\theta_1 = O(\alpha)$.
We now combine (\ref{expansionOmega}) with
\beq \label{Omegaat0}
\Omega(0) = -  {\psi(0) \over c [ U'_s(0) \psi(0) + c \psi'(0) ]} 
\eeq
and obtain
\beln \label{valueOmega}
\Omega(0)^{-1} & = \alpha (U_+ - c)^2 - \alpha^2 (U_+ - c)^4 \Omega_0(0,c) + O(\alpha^3) 
 \cr & =  - c U_s'(0) - c^2 {\psi_{-,\alpha}'(0) \over \psi_{-,\alpha}(0)},
\eeln
which leads to (\ref{psippsi}) and ends the proof of Theorem \ref{theo3}.


\subsection{Construction on $[-1,+1]$ \label{interval}} 


Let $U_s(y)$ be a $C^\infty$ function on $[-1,+1]$, which is even, and such that $U_s(\pm 1) = 0$.
Let $\psi(y,c)$ be the even solution of Rayleigh equation. In this section we compute $\psi'(1,c) / \psi(1,c)$.

We introduce 
$$
\omega(y) = {1 \over \Omega(y)} = (U_s(y) - c) { U_s'(y) \psi(y) - (U_s(y) - c) \psi'(y) \over \psi(y)} .
$$
Then a direct computation shows that
\beq \label{equationomega}
\omega'(y) = - \alpha^2 Y(y) + {\omega^2(y) \over Y(y)}
\eeq
where
$$
Y(y) = (U_s(y) - c)^2 .
$$
Moreover, $\omega(0) = 0$ since $\psi$ is even.
We look for an expansion of $\omega(y)$ of the form
$$
\omega(y) = \omega_2(y) \alpha^2 + \omega_4(y) \alpha^4 + O( \alpha^6) .
$$
This leads to
$$
\omega_2'(y) = Y(y)
$$
and to
$$
\omega_4'(y) = {\omega_2^2(y) \over Y(y)}.
$$
The solutions which vanish at $y = 0$ are
$$
\omega_2(y) = \int_y^0 Y(z) \, dz
$$
and
$$
\omega_4(y) = - \int_y^0 {\omega_2^2(z) \over Y(z) } dz.
$$
This leads to
\beq \label{expansionomega1}
\omega(1) = - \alpha^2  \int_0^1 (U_s(z) - c)^2 \, dz + \alpha^4  \int_0^1 {\omega_2^2(z) \over (U_s(z) - c)^2} \, dz + O(\alpha^6).
\eeq


\section{Green function \label{secGreen}}


We now construct the Green function $G(x,y)$ of Rayleigh equation, which by definition satisfies
$$
Ray(G(x,y)) = \delta_x
$$
with boundary conditions $G(x,0) = 0$ and $G(x,y) \to 0$ as $y \to + \infty$.
Let $\psi_\pm(y)$ be two solutions of Rayleigh equation with non zero Wronskian, such that $\psi_-(y) \to 0$ as $y \to + \infty$.
We will construct the Green function in two steps. First we construct an "interior" Green function $G^{int}(x,y)$ which satisfies Rayleigh equation
and $G^{int}(x,y) \to 0$ as $y \to + \infty$. Next we add a "boundary" Green function $G^b(x,y)$ to recover $G(x,0) = 0$.

For the "interior" Green function $G^{int}(x,y)$, we may choose
\beq \label{GreenRayint}
G^{int}(x,y) = {1 \over U_s(x) - c} {1 \over W[\psi_-,\psi_+]} \Bigl\{
\begin{array}{c}
\psi_+(x) \psi_-(y) \qquad \hbox{if} \qquad x < y, \cr
\psi_-(x) \psi_+(y) \qquad \hbox{if} \qquad x > y. 
\end{array}
\eeq
Then, by construction, $\psi^{int}$, defined by
$$
\psi^{int}(y) = \int_{\rit^+} G^{int}(x,y) f(x) \, dx 
$$
satisfies
$$
Ray(\psi^{int}) = f,
$$
together with $\psi^{int}(y) \to 0$ as $y \to + \infty$. Note that $G^{int}(x,y)$ is well defined for any right-hand side $f$.

We now introduce $G^b$ such that $G^{int}(x,0) + G^b(x,0) = 0$. We can choose
\beq \label{GreenRayb}
G^b(x,y) = {1 \over W[\psi_-,\psi_+]} {\psi_+(0) \over \psi_-(0)} \psi_-(x) \psi_-(y).
\eeq
We then define 
$$
\psi^b(y) = - {1 \over W[\psi_-,\psi_+]} {\psi_+(0) \over \psi_-(0)} \Bigl[ \int_{\rit_+} \psi_-(x) f(x) \, dx \Bigr] \psi_-(y).
$$
By construction,
$$
Ray(\psi^{int} + \psi^b) = f
$$
and $\psi^{int} + \psi^b$ satisfies the boundary conditions at $y = 0$ and $y = + \infty$.

Note that the boundary Green function is singular when $\psi_-(0) = 0$, namely near eigenvalues of Rayleigh operator.


\subsubsection*{Acknoledgements}


The authors would like to thank N. Masmoudi and W. Zhao for many discussions on Rayleigh equation, in particular during an invitation
at NYU Abu Dhabi in $2024$.
 D. Bian is supported by NSFC under the contract 12271032.



\end{document}